\documentclass{cmslatex}
\usepackage{latexsym, amssymb, enumerate, amsmath}
\usepackage[dvips]{epsfig}  

\renewcommand{\theequation}{\arabic{section}.\arabic{equation}}

\newtheorem{thm}{Theorem}[section]

\newtheorem{remark}{Remark}[section]  

\newtheorem{rem}[thm]{Remark}

\begin{document}

\title{Kinetic derivation of a Hamilton-Jacobi  traffic flow model
\thanks{
}}
\author{Raul Borsche 
\thanks {Fachbereich Mathematik, Technische Universit\"at
  Kaiserslautern, Germany, (borsche@mathematik.uni-kl.de).}
\and  Marc Kimathi\thanks {Fachbereich Mathematik, Technische Universit\"at
  Kaiserslautern, Germany, (kimathi@mathematik.uni-kl.de).}
 \and Axel Klar\thanks {Fachbereich Mathematik, Technische Universit\"at
  Kaiserslautern, Germany, (klar@mathematik.uni-kl.de).} }
  



\pagestyle{myheadings} \markboth{Hamilon-Jacobi traffic model}{R. Borsche, M. Kimathi, A. Klar}\maketitle

\begin{abstract}
Kinetic  models for vehicular traffic are reviewed and considered from the point of view of deriving  macroscopic equations.
A derivation of the associated
macroscopic traffic flow equations leads to different types of equations:  in certain situations modified Aw-Rascle equations are obtained. On the other hand,
for several choices of kinetic parameters new Hamilton-Jacobi type traffic equations are found. Associated microscopic models are discussed and 
numerical experiments are presented discussing several situations
 for highway traffic and comparing the different models.
 \end{abstract}


\section{Introduction}
\label{Introduction}

Macroscopic models for vehicular traffic have been first introduced by  Lighthill and Whitham \cite{Whi74}. These models are based on  the continuity
equation for the density $\rho$ closing the equation by an equilibrium
assumption on the mean velocity $u$, where  $u$ is approximated by
a uniquely determined 
 equilibrium value, \cite{Whi74}. An additional momentum equation for $u$
has been introduced by Payne and Whitham in \cite{Pay79,Whi74} in
analogy to fluid dynamics. To avoid certain inconsistencies, like wrong way traffic, of  models such as the
Payne/Whitham model a new macroscopic model has been introduced by Aw and
Rascle \cite{AR98}, see also  \cite{aw00} or
 \cite{Gre00}. These models have been subsequently improved, for example, in \cite{Deg,Ber}.

Kinetic equations for vehicular traffic can be found, for example, in
\cite{PH71,PF75,Nel95,KW97,IKM03}. Procedures to derive macroscopic traffic
equations including the Aw/Rascle model from underlying kinetic models
have been performed in different ways by several authors, see, for
example, \cite{Hel95B} and \cite{KW00}. These procedures are developed
in analogy to the transition from the kinetic theory of gases to
continuum gas dynamics.  

In the present paper these derivations are reviewed. A closer analysis shows that Aw-Rascle type traffic equations can be derived from kinetic problems
for certain choices of kinetic parameters. For other choices, however, new equations with  Hamilton-Jacobi  terms are derived.

The paper is arranged in the following way: In Section \ref{The
Kinetic Model} different reduced  kinetic models are presented.   Section \ref{Derivation of Macroscopic
Models} contains the derivation of the macroscopic models mentioned above.
Section \ref{Associated} contains the associated microscopic traffic flow models.
  Finally, in Section
\ref{Numerical Investigations} numerical results are given comparing the derived  macroscopic equations  for several nonhomogeneous traffic
flow situation.

\section{Kinetic  Models}
\label{The Kinetic  Model}

The kinetic models presented in this section  are based on  work in \cite{KW981,IKM03}
and describe highway traffic in a cumulative way averaging over all
lanes. These models are given by integro-differential and Fokker-Planck type equations respectively. 
In particular, the Fokker-Planck type models are changed allowing only for densities below  a maximal density and 
for a better comparison of the models.

\subsection{Correlations and the reduced density}

The basic quantity in a kinetic approach is the single car
distribution $f (x,v)$ describing the (number) density of cars at $x$ with
velocity $v$. The total density $\rho$ on the highway is
defined by
\begin{eqnarray*}
 \rho (x) \;=\; \int_0^w  f(x,v) dv,
\end{eqnarray*}
where $w$ denotes the maximal velocity.  Let $F(x,v)$ denote the
probability distribution in $v$ of cars at $x$, i.e. $f(x,v) = \rho(x)
F(x,v)$.  Then , the  mean velocity is
\begin{eqnarray*}
 u (x) \;=\; \int_0^w v  F(x,v) dv.
\end{eqnarray*}

An important role is played by the distribution $f^{(2)} (x,v,h,v_+)$
of pairs of cars being at the spatial point $x$ with velocity $v$ and
leading cars at $x+h$ with velocity $v_+$.  This distribution function
has to be approximated by the one-vehicle distribution function
$f (x,v)$. Usually, a  chaos assumption is used,
\begin{eqnarray*}
f^{(2)} (x,v,h,v_+) = q(h,v;f)\, f(x,v)\, F(x+h,v_+),
\end{eqnarray*}
compare Nelson \cite{Nel95}.  For a vehicle with velocity $v$ the
function $q(h,v;f)$ denotes the distribution of leading vehicles
with distance $h$ under the assumption that the velocities of the
vehicles are distributed according to the distribution function $f$.

Moreover, we introduce thresholds for braking ($H_{B}$)
and acceleration ($H_{A}$):
\begin{eqnarray*}
H_X = H_{X} (v) \;=\; H_0 + v T_{X}, \quad X = B,A.
\end{eqnarray*}
$T_{B} < T_{A}$ are
reaction times. $H_0$ denotes the
minimal distance between the vehicles.   
For simplicity we choose $H_A$ and $H_B$ in the following as constants.

The distribution of leading vehicles $q(h,v;f)$ is prescribed a
priori.  The main properties, which $q(h,v;f)$ has to fulfill are
positivity, 
\begin{eqnarray*}
\int_0^\infty q(h,v;f)dh \;=\; 1, 
\end{eqnarray*}
and
\begin{eqnarray}
\label{meanvalue}
\int_0^w \int_0^{\infty} h q(h,v;f) dh\, F(v) dv  \;=\;  \frac{1}{\rho}.
\end{eqnarray}
Equation (\ref{meanvalue})
 means that the average headway of the cars is $1/\rho$.
Here, the leading vehicles are assumed to be distributed in an uncorrelated
way with a minimal distance $H_B$ from the car under
consideration, see Nelson \cite{Nel95}:
\begin{eqnarray*}
q(h,v;f) \;=\;q(h;\rho) =  \tilde{\rho}\, e^{-\tilde{\rho}(h-H_B)}\,\chi_{[H_B,\infty)}(h).
\end{eqnarray*}
The reduced density $\tilde\rho$ has to be defined in such a way, that
(\ref{meanvalue}) is fulfilled. One obtains
\begin{eqnarray}
\label{admissible}
\tilde\rho\;=\;\frac{\rho}{1-\rho\,\int_0^w H_B F (v)  dv}
\;=\;\frac{\rho}{1-\rho\, H_B}.
\end{eqnarray}

\begin{remark} 
The reduced density $\tilde\rho$ must be positive, i.e.
\begin{eqnarray*}
\rho\; \;<\; \frac{1}{H_B}
\end{eqnarray*}
\end{remark}

We note that 
\begin{eqnarray*}
q(H_A;\rho) 
\;=\; \tilde\rho\,e^{-\tilde\rho(H_A-H_B)}
\end{eqnarray*}
and
\begin{eqnarray*}
q(H_B;\rho) \;
\;=\; \tilde\rho\;.
\end{eqnarray*}

Moreover, from phenomenological considerations the probability of braking can be derived as
$$
P_B = 1-(1-\rho H_B) e^{- \tilde \rho H_B},
$$
see \cite{GKMW03}.
These basic considerations can be used to develop different kinetic models.

\subsection{Models based on Integro-differential equations}

A first  kinetic  model is derived using  classical Boltzmann arguments. It  is given by the following evolution equation
for the distribution function $f$, see \cite{KW00,GKMW03}:

\begin{eqnarray}
\label{kinetic2}
\partial_t  f + v \partial_x  f
&=&  C^+ (f)\\
&=&
 \left[  q_B P_B (G^+_B -L^+_B ) (f) + q_A (G^+_{A} -L^+_{A}) (f)  
\right]  \nonumber
\end{eqnarray}    
with
\begin{eqnarray*}
G_{B}^+ (f)
&=&
\int \int_{\hat{v} > \hat{v}_+}
\vert \hat{v} -\hat{v}_+ \vert
\sigma_{B} (v ; \hat{v},\hat{v}_+) 
 f(x,\hat{v})
 F(x+H_B, \hat{v}_+)
d \hat{v} d\hat{v}_+\\
L_{B}^+ (f)
&=&
\int_{\hat{v}_+<v} 
\vert v -\hat{v}_+ \vert
f(x,v)
 F(x+H_B, \hat{v}_+) d\hat{v}_+\\
G_{A}^+ (f)&=&
\int \int_{\hat{v} < \hat{v}_+}
 \vert \hat{v} -\hat{v}_+ \vert
\sigma_{A} (v ; \hat{v}, \hat{v}_+) 
 f(x,\hat{v})
 F(x+H_A, \hat{v}_+)
d \hat{v} d\hat{v}_+\\
L_{A}^+ (f)&=&
\int_{\hat{v}_+>v}  
\vert v -\hat{v}_+ \vert 
 f(x,v)
 F(x+H_A, \hat{v}_+)  d\hat{v}_+
 \end{eqnarray*}  

$G_B, L_B$ stand for gain and loss terms resulting from braking interactions, $G_A, L_A$ result from accelerating interactions.
 Reaching the braking line the vehicle
brakes, such that the new velocity $v$ is distributed with a
distribution function $\sigma_B$ depending on the old velocities
$\hat{v},\hat{v}_+$.
For acceleration, the new velocity is distributed according to $\sigma_A$.

\begin{remark}
In \cite{GKMW03} additionally a relaxation term is
introduced, describing a random behaviour of the drivers. It is given
by
\begin{eqnarray*}
 G_{S} (f) -  L_{S} (f) \;=\;
\nu  ( \int_{0}^w
\sigma_{S} (v , \hat{v})
f(x, \hat{v})
d \hat{v}
-   f(v) )  .
\end{eqnarray*}    
This term is necessary as long  as one is interested in a more detailed investigation of the stationary solutions of the kinetic model
and the resulting fundamental diagrams.
However, in the present investigation we aim at deriving different macroscopic equations without relaxation terms on the right hand side.
For such a derivation it is sufficient to  consider the simplified version above.  
For further remarks on this Boltzmann/Enskog approach to traffic flow
modelling see \cite{KW00}.  
\end{remark}

\subsection*{Example 1}

For the  probability distributions $\sigma_A, \sigma_B$
we  choose the following simple expressions, see \cite{GKMW03}:
\begin{eqnarray}
\label{1s}
\sigma_{B} (v , \hat{v},\hat{v}_+) \;=\; \frac{1}{\hat{v} - \hat{v}_+}
\,\chi_{[\hat{v}_+, \hat{v}]} (v)
\end{eqnarray}
and
\begin{eqnarray}
\label{2s}
\sigma_{A} (v , \hat{v} , \hat{v}_+) \;=\; \frac{1}{\hat{v}_+ -  \hat{v}}
\,\chi_{[\hat{v},  \hat{v}_+]} (v).
\end{eqnarray}
This means we have an equidistribution of the new velocities
between the velocity of the car and the velocity of its leading
car.

\subsection*{Example 2}

Another possible choice is, see \cite{KW00}
\begin{eqnarray*}
&&\sigma_{B} (v , \hat{v}) = \frac{1}{\hat{v} (1-\beta)}
\chi_{[\beta  \hat{v},  \hat{v}]} (v)
\end{eqnarray*}
and 
\begin{eqnarray*}
\sigma_{A} (v , \hat{v}) = \frac{1}{\mbox{min}(w,\alpha \hat{v}) -  \hat{v}}
\chi_{[\hat{v},  \mbox{min}(w,\alpha \hat{v})]} (v).
\end{eqnarray*}

\subsection{Models based on Vlasov-Fokker-Planck equations}

In \cite{IKM03} a kinetic model based on a VlasovFokker-Planck approach has been developed:
\begin{equation}
\label{fp}
   \partial_t f+v\partial_x f  = C^+ (f) = - \partial_v\big(B[f]f\big) .
\end{equation}

Here, $f$ stands again for a traffic distribution function. We denote by  $\rho,u$ again the macroscopic density
and speed associated with $f.$

To define the braking and acceleration behaviour of drivers in response to 
traffic situations, we use  the following braking/acceleration forces
as functions of the traffic conditions. Slightly changing the approach in \cite{IKM03}, i.e. adding the parameters
$q_B$ and $q_A$,  we consider 
\begin{equation} \begin{split}
   B[f](t,x,v) \;=\; \begin{cases}
                            - q_B P_B c_\eta \vert v-u^B \vert^\eta\, & \quad v\,>\,u^B \\[2mm]
		           q_A c_\eta \vert u^A-v\vert ^\eta & \quad v\,\leq\,u^B \text{ and } v\,\leq\,u^A \\[2mm]
                           0 & \quad \text{else}
	               \end{cases}
\end{split} \label{power_term} \end{equation}
Again we look at two examples, i.e.   $\eta =1$ and $\eta =2$. Here  $c_\eta = v_{ref}$ with $v_{ref}$ a reference velocity if $\eta =1$ and $c_\eta$ dimensionless if $\eta =2$ 
and 
\begin{align}\begin{split}
   \rho^X \;=\; \rho(x+H_X,t),\quad\quad & u^X \;=\; u(x+H_X,t) 
\end{split}\end{align}
for $X=A,B$. 

\begin{remark}
Similar to the case of the integro-differential equation, we use for the present investigation  a simplified version of the kinetic model, see also \cite{IH}. 
In the original version of the model in \cite{IKM03} a diffusion term 
\begin{eqnarray*}
 \partial_v (D[f] \partial_v f ) 
 \end{eqnarray*}
 with
\begin{eqnarray*} 
  D[f](\rho,u,v) = \begin{bmatrix} \sigma(\rho^B,u^B) |v-u^B|^\gamma &
 v>u^B \\
\sigma(\rho^A,u^A) |v-u^A|^\gamma &
 \mbox{ else }  \end{bmatrix}
\end{eqnarray*}
with $\gamma\ge 1$ has been added to the right hand side of the above equation.  Details of the  function $\sigma(\rho,u)$ can be found in 
 reference~\cite{IKM03}.
For the presentation here we neglect this diffusion term. It is however  necessary to obtain smooth homogeneous solutions.
\end{remark}

\section{Derivation of  Macroscopic  Models}
\label{Derivation of Macroscopic Models}

In this section macroscopic
equations for density and mean velocity are derived
following the procedure in \cite{KW00}. Among these equations are new Hamilton-Jacobi type traffic equations which have not been discussed up to now in literature.
This section shows that the resulting equations do not depend on the the different kinetic models used, but rather on the type of interaction terms. Using simplified closure relations explicit results are obtained compared to the numerical closures in \cite{KW00}. However, the resulting equations are still more detailed than the usually used macroscopic models.

\subsection{Balance Equations}
\label{Balance Equations}

Multiplying the inhomogeneous kinetic equation (\ref{kinetic2}) or  (\ref{fp}) with
$1$ and $v$ and integrating it with respect to $v$ one obtains the
following set of balance equations:
\begin{eqnarray}
\label{balance}
\partial_t \rho + \partial_x (\rho u) & = & 0\\
\partial_t (\rho u) + \partial_x ( P + \rho u^2)
+ E
& = & 0\nonumber
\end{eqnarray}
with the 'traffic pressure'
\begin{eqnarray}
\label{pressure}
P \;=\;\int_0^w  (v-u)^2 f dv,
\end{eqnarray}
and the  flux term
\begin{eqnarray}
\label{ens}
E \;=\; - \int_0^w v  C^+(f) (x,v,t) dv.
\end{eqnarray}

To obtain closed equations for $\rho$ and $u$ one has to specify the
dependence of $P$ and $E$ on $\rho$ and $u$.

\subsection{Closure  and resulting macroscopic equations}
\label{Closure Relations}

To approximate the distribution function we use the simplest possible one node quadrature ansatz disregarding fluctuations in the distribution function. That means, we use
 $f (v) \sim  \rho \delta_u (v)  $ for
the distribution function in (\ref{pressure}) and (\ref{ens})  to approximate the true distribution $f$ and
to close the equations, compare  \cite{IH} for such an ansatz in the traffic case or \cite{CDP} for a similar procedure for interacting particle systems.  
Using this Ansatz, one obviously neglects the variance of the distribution function. However, the main features of the resulting 
macroscopic equation are preserved. We obtain for  the traffic
pressure 
$$P \sim  0.$$

We are left with the Enskog term $E$. It is approximated by considering 
expression (\ref{ens}) for $E$  and substituting
the closure  for $f$.  One obtains different expressions depending on the kinetic model under consideration.

\subsubsection{Integro-differential equations}
In the case of integro-differential equations one obtains
\begin{eqnarray*}
E \;=\; E_B(f)+E_A(f)
\end{eqnarray*}
with
\begin{eqnarray*}
E_B(f)
& = & - q_B P_B \int \int_{\hat{v} > \hat{v}_+}
\vert \hat{v}-\hat{v}_+ \vert\\
&&  f(x, \hat{v})  F(x+H_B , \hat{v}_+)
[\int_0^w v \sigma_B(v, \hat{v},\hat{v}_+) dv - \hat{v}]
 d\hat{v}_+ d\hat{v} \\
\end{eqnarray*}
and
\begin{eqnarray*}
E_A(f) &=&-  q_A 
\int \int_{\hat{v}  <  \hat{v}_+}
\vert \hat{v}-\hat{v}_+ \vert \\
&&f(x,\hat{v}) F(x +H_A , \hat{v}_+)
[\int_0^w v \sigma_A (v, \hat{v},\hat{v}_+) dv - \hat{v}]
d\hat{v}_+ d\hat{v}.
\end{eqnarray*}

Using now $$ F (x,v)  = \delta _{u(x)}  (v)$$  gives for $ u >  u^B$ approximately:
\begin{eqnarray*}
E_B \sim - q_B P_B \rho \vert u -  u^B\vert  [\int_0^w v \sigma_B(v, u , u^B) dv - u ]     
\end{eqnarray*}
and $0$ otherwise. Approximating $u^B -  u$ by $H_B \partial_x u$  this is approximated for $ \partial_x u < 0$ by 
\begin{eqnarray*} q_B P_B \rho H_B \partial_x u   [\int_0^w v \sigma_B(v, u , u^B) dv - u ] .  
\end{eqnarray*}
The acceleration term gives 
\begin{eqnarray*}
E_A  \;\sim\; - q_A \rho \vert u -  u^A \vert  [\int_0^w v \sigma_A(v, u , u^A) dv - u ]  
\end{eqnarray*}
for $ u <  u^A$ and $0$ otherwise. Therefore one obtains for $\partial_x u > 0$ the approximation
\begin{eqnarray*}
   - q_A  \rho H_A \partial_x u  [\int_0^w v \sigma_A(v, u , u^A) dv - u ]  .
\end{eqnarray*}
The final result depends on the interaction model.
Example 1 gives
$$E = 
 \left\{\begin{array}{ll}
             E_B \sim - q_B P_B \rho H_B^2 \partial_x u \vert  \partial_x u \vert ,& \partial_x u<0 \\
             E_A \;\sim\;  - q_A  \rho H_A^2 \partial_x u   \vert \partial_x u \vert, & \partial_x u > 0.
         \end{array}
         \right.
$$
Example 2 gives
$$E =
    \left\{\begin{array}{ll}
             E_B \sim - q_B P_B \rho H_B    \frac{1-\beta}{2} u\partial_x u,& \partial_x u<0 \\
             E_A \;\sim\;  - q_A \rho  H_A \frac{min(\alpha u, w) -u}{2} \partial_x u,& \partial_x u > 0.
         \end{array}
         \right.
$$
\subsubsection{Vlasov-Fokker-Planck equations}
Similar results are obtained for the Vlasov-Fokker-Planck equations.
Computing
\begin{eqnarray*}
E  =   \int v  \partial_v\big(B[f]f\big)   dv =  - \int \big(B[f]f\big)   dv \end{eqnarray*}
one obtains for $ u >  u^B$
\begin{eqnarray*}
E \sim  c_\eta q_B P_B \rho \vert u - u^B \vert^\eta    
\end{eqnarray*}
and for $ u <  u^B$ and $u < u^A$
\begin{eqnarray*}
E \;\sim\;  -  c_\eta q_A  \rho  \vert u - u^A \vert^\eta   \end{eqnarray*}
and $0$ else.
This gives for $\eta=1 $
$$E \sim 
 \left\{\begin{array}{ll}
             - v_{ref} q_B P_B \rho H_B \partial_x u, & \partial_x u<0 \\
             -  v_{ref} q_A  \rho H_A   \partial_x u,& \partial_x u > 0.
         \end{array}
         \right.
$$
For $\eta=2$ we have 
$$E \sim 
 \left\{\begin{array}{ll}
             - c_\eta q_B P_B \rho H_B^2 \vert \partial_x u  \vert  \partial_x u ,& \partial_x u<0 \\
             - c_\eta q_A  \rho H_A^2 \vert \partial_x u  \vert  \partial_x u, & \partial_x u > 0.
         \end{array}
         \right.
$$
\begin{rem}
In both cases, depending on the interaction law, either a linear dependence on $\partial_x u $ or a nonlinear functional dependence is observed.
\end{rem}
\subsection{Macroscopic equations}

Altogether, one obtains macroscopic equations either of the form
\begin{eqnarray}
\label{macro}
\partial_t \rho + \partial_x (\rho u) &=& 0\\
\partial_t (\rho u) + \partial_x ( \rho u^2)
-    \rho a(\rho,u) \partial_x u
&=& 0 \nonumber
\end{eqnarray} 
or of the form
\begin{eqnarray}
\label{HJ}
\partial_t \rho + \partial_x (\rho u) &=& 0\\
\partial_t (\rho u) + \partial_x ( \rho u^2)
-   \rho b(\rho,u) \vert \partial_x u \vert \partial_x u 
&=& 0 ,\nonumber
\end{eqnarray} 
where the  coefficients are given by 
$$a(\rho,u) = 
    \left\{\begin{array}{ll}
             \frac{H_B P_B}{\frac{1}{\rho}-  H_B}  f_B(u) & \partial_x u<0 \\
             \frac{H_A}{\frac{1}{\rho}-  H_B} \exp(- \tilde \rho (H_A- H_B)) f_A(u) & \partial_x u > 0
         \end{array}
         \right.
$$
$$b(\rho,u) = 
 \left\{\begin{array}{ll}
           \frac{H_B^2 P_B}{\frac{1}{\rho}-  H_B}  & \partial_x u<0 \\
            \frac{H_A^2}{\frac{1}{\rho}-  H_B}  \exp(- \tilde \rho (H_A- H_B)) & \partial_x u > 0
         \end{array}
         \right.
$$
with  suitable functions $f_A, f_B $. We note that $a(\rho,u), b (\rho,u) >0$.
Looking at these equations one observes that equation (\ref{macro}) 
is a  Rascle-type equation with microscopically justified coefficients which include braking and acceleration threshold.
On the other hand,  equation (\ref{HJ})
is an  equation with Hamilton-Jacobi terms, which has, to the knowledge of the authors, not been discussed in the literature.  Vehicles described by (\ref{HJ}) will brake stronger or accelerate faster, the steeper the gradient in velocity is ahead of them. 

If we simplify further, choosing $H_A = H_B =H$ and $q_A = q_B = \tilde\rho$, $P_B =1$ and approximating $f_A,f_B$ by $v_{ref}$
one obtains the coefficients  
\begin{eqnarray}
\label{a} a(\rho) = \frac{H v_{ref}}{\frac{1}{\rho}-  H} =  \frac{v_{ref}}{\frac{1}{\rho H }-  1} \\
\label{b}
b(\rho) =  \frac{H^2}{\frac{1}{\rho}-  H}=  \frac{H}{\frac{1}{\rho H }-  1}.
\end{eqnarray}

\begin{remark}
Equation (\ref{macro}) with the coefficient (\ref{a})  is similar to the modified Rascle equation discussed together with its limits in   \cite{Deg}.
From the kinetic point of view these equations are strongly simplified. In particular, they treat the braking and acceleration interaction in the same way, which is clearly not physical. However, they still contain the essential features of traffic flow, see 
\cite{Deg}.
\end{remark}

\begin{remark}
The kind of equation one obtains does not depend on the fact whether an integro-differential equation model or a Fokker-Planck type model is used, but rather on the fact 
which interaction rule is chosen.
\end{remark}
\begin{remark}
We note that traffic equations with different Hamilton-Jacobi terms have also been discussed in \cite{IH}. 
\end{remark}

\begin{remark}
The two results obtained here could be also merged into a third equation by using 
\begin{eqnarray}
\label{merge}
\partial_t \rho + \partial_x (\rho u) &=& 0\\
\partial_t (\rho u) + \partial_x ( \rho u^2)
-    \rho b (\rho) c(\vert \partial_x u \vert) \partial_x u
&=& 0 \nonumber
\end{eqnarray} 
with
\begin{eqnarray*}
 c(  \vert \partial_x u \vert ) =  \mbox{min} \;  \{  \vert \partial_x u \vert, C \}
\end{eqnarray*}
with a constant $C$.
This would limit   the braking force.
\end{remark}

\section{Associated microscopic car-following models}
\label{Associated}
Equation (\ref{macro}) with coefficient (\ref{a}) can be  derived from microscopic models of the form
\begin{eqnarray*}
\dot{x}_i &=& v_i \\
\dot{v}_i &= &\frac{H v_{ref}}{x_{i+1} - x_i }  \frac{v_{i+1}- v_i}{x_{i+1} - x_i - H}.
\end{eqnarray*}
This can be easily seen by the following procedure, compare \cite{AKMR03}.
Set \begin{eqnarray*}
l_i &=& x_{i+1} - x_{i},
\end{eqnarray*}
then  the microscopic  equations are
\begin{eqnarray*}
\dot{x}_i &=& v_i\\
\dot{v}_i &=& \frac{H v_{ref}}{l_i }  \frac{v_{i+1}- v_i}{l_i - H}.
\nonumber
\end{eqnarray*}
The local 
(normalized) density around vehicle i and its inverse 
the local (normalized) specific volume are respectively 
defined by
\begin{eqnarray*}
\rho_i = \frac{H}{l_i} \; \mbox{and} \;\; 
\tau_i = \frac{1}{ \rho_i} = \frac{l_i}{H}.
\end{eqnarray*}
One obtains
the  microscopic model
\begin{eqnarray}
\label{ode2}
\dot{x}_i &=& v_i \;,\\
\dot{v}_i &=& \frac{v_{ref}}{\tau_i}  \frac{1}{H} \; \frac{(v_{i+1} -v_i) }
{\tau_{i}-1}
 . \nonumber
\end{eqnarray}
We have
\begin{eqnarray*}
\dot{l}_i = v_{i+1}-v_i
\;\; \mbox{or} \;\; 
\dot{\tau}_i = \frac{1}{H} \; (v_{i+1}-v_i).
\end{eqnarray*}  
One considers the coordinate $X= \int^x \rho(y,t) dy$ describing
the 
total space occupied by cars up to point $x$.
 Approximating $(v_{i+1} - v_i )/H$ by $\partial_X u$ yields the Lagrangian form of the macroscopic equations, i.e. the equivalent of the p-system in gas dynamcis
\begin{eqnarray}
\label{pde}
\partial_T \tau - \partial_X u &=& 0 \; \\
\partial_T u - \frac{a (\tau) }{\tau} \partial_X u &=& 0\; ,
\nonumber
\end{eqnarray}
where
\begin{equation}\label{ptau}
a (\tau)=
 \frac{v_{ref}}{\tau-1}.
\end{equation}

We change the Lagrangian
``mass'' coordinates $(X,T)$ into Eulerian coordinates  $(x,t)$ with
$$\partial_x X= \rho , \;\partial_t X = - \rho v , \;
T = t$$
or
$$\partial_X x= \rho^{-1} = \tau,  \;
\partial_T x =  v.$$
The macroscopic 
system in Eulerian coordinates is then 
\begin{eqnarray}
\label{cons}
\partial_t \rho + \partial_x (\rho u ) &=&0 ,\\
 \partial_t (\rho u)  + \partial_x ( \rho u^2 ) - \rho a(\rho) \partial_x u  &=& 0 \nonumber
\end{eqnarray}
with
\begin{equation}\label{prho}
a (\rho)= v_{ref}
  \left( \frac{1}{\rho } -1 \right)^{-1} .
\end{equation}
This means we obtain again the  equations (\ref{macro}) and  (\ref{a}) taking into account that in the kinetic derivation $\rho$ is the number density. That means 
the quantity $\rho H$ in the kinetic part is  equivalent to the normalized density considered in this section.

\begin{remark}
We note that the above statement is equivalent to considering the kinetic equations for the rescaled distribution functions $f^\prime = f H$. This leads, for example, to a Vlasov equation where the braking and acceleration term in (\ref{fp}) is multiplied by $\frac{1}{H}$.

\end{remark}

\begin{remark}
For numerical simulations of the microscopic system and comparison with the macroscopic equation the quantity $H$ is  chosen such that the total space $\int_0^L \rho(x) dx $ occupied by the cars  is equal to  $H N$, where $L$ is the total length of the region under consideration and $N$ is the total number of vehicles.

\end{remark}

Using the same procedure one obtains the  microscopic model associated to equation (\ref{HJ}) and (\ref{b}). It is given by 
\begin{eqnarray*}
\dot{x}_i &=& v_i \\
\dot{v}_i &= &\frac{ H}{(x_{i+1} - x_i)^2 }. \frac{ \vert v_{i+1}- v_i \vert (v_{i+1}- v_i)}{x_{i+1} - x_i - H}.
\end{eqnarray*}
The latter equations are similar to microscopic traffic equations originally stated by Wiedemann and Leutzbach \cite{WL}.

\section{Numerical Investigations}
\label{Numerical Investigations}
 
In this section we investigate the macroscopic equations numerically. In particular, the Hamilton-Jacobi type equations equation (\ref{HJ}) are compared to the Aw-Rascle type equations (\ref{macro}) .

\subsection{Numerical methods}

We choose a numerical method suited for the hyperbolic equation in non-conservative  form (\ref{macro}) as well as for the Hamilton-Jacobi term in (\ref{HJ}).
A suitable choice is given e.g. by second order central scheme developed in \cite{BL03}.
For completeness we state an extended version of the scheme as used in our numerical computations. 
To start with, the above equations are written in the form
\begin{eqnarray}
\label{HJnum}
\partial_t \phi + H(\phi,  \phi_x) &=&0
\end{eqnarray}
 with
 $$\phi = 
 \left(\begin{array}{ll}
           \rho \\
          \rho u 
         \end{array}
         \right).
$$ 
For equations (\ref{macro}) we have 
 $$H(\phi, \phi_x) = 
 \left(\begin{array}{ll}
           \rho u_x  + u \rho_x\\
          (\rho u^2)_x - \rho a(\rho) u_x)
         \end{array}
         \right)
$$ 
and for equations   (\ref{HJ})
$$H(\phi, \phi_x) = 
 \left(\begin{array}{ll}
           \rho u_x  + u \rho_x\\
          (\rho u^2)_x - \rho b(\rho) \vert u_x \vert u_x)
         \end{array}
         \right).
$$ 
For the numerical scheme a grid of equally spaced points $x_i\ i=1,\dots,N$, with $\Delta x=x_i-x_{i-1}$ is given.
In the following we consider the explicit time step from $t_m$ to $t_{m+1}=t_m+\Delta t$.
The  aim is to construct a second order scheme for the above 1-D equations.
A detailed derivation can be found in \cite{BL03}. 

Based on piecewise quadratic interpolations 
one obtains the following expression for the iterate $\phi_i^{m} $ approximating $\phi(x_i,t_m)$;
\begin{eqnarray}
\label{scheme}
\phi_i^{m+1} = \phi_{i-\frac{1}{2}}^{m+1} + \frac{1}{2} (\Delta \phi)^{m+1}_i - \frac{1}{8} \mathcal{D} (\Delta \phi)^{m+1}_i
\end{eqnarray}
with the second order approximation of the equation
\begin{eqnarray*}
\phi^{m+1}_{i-\frac{1}{2}} = \phi_{i-\frac{1}{2}}^{m} - \Delta t  H(\phi_{i-\frac{1}{2}}^{m+\frac{1}{2}}, (\phi_x)_{i-\frac{1}{2}}^{m+\frac{1}{2}})
\end{eqnarray*}
where
\begin{eqnarray*}
\phi_{i-\frac{1}{2}}^{m+\frac{1}{2}}&=&\phi_{i-\frac{1}{2}}^{m}- \frac{\Delta t }{2} H(\phi_{i-\frac{1}{2}}^{m},(\phi_x)_{i-\frac{1}{2}}^{m})\\
 (\phi_x)_{i-\frac{1}{2}}^{m+\frac{1}{2}}&=&(\phi_x)_{i-\frac{1}{2}}^{m} - \frac{\Delta t }{2} 
 [\frac{\partial H}{\partial \phi} (\phi_{i-\frac{1}{2}}^{m},(\phi_x)_{i-\frac{1}{2}}^{m})  (\phi_x)_{i-\frac{1}{2}}^{m} \\
&& + \frac{\partial H}{\partial \phi_x} ( \phi_{i-\frac{1}{2}}^{m},(\phi_x)_{i-\frac{1}{2}}^{m})
 \frac{\mathcal{D} (\Delta \phi)^{m}_{i-\frac{1}{2}} }{(\Delta x )^2} ].
  \end{eqnarray*}
In these expressions the following definitions are obtained from Taylor expansions:
\begin{eqnarray*}
\phi^{m}_{i\pm \frac{1}{2}} &=& \phi_{i}^{m} \pm  \frac{1}{2}  (\Delta \phi)_{i \pm \frac{1}{2}}^{m}- \frac{1}{8}  \mathcal{D} (\Delta \phi)^{m}_{i\pm\frac{1}{2} }\\
(\phi_x)_{i-\frac{1}{2}}^{m}  &=& \frac{(\Delta \phi)_{i-\frac{1}{2}}^{m}}{\Delta x }
\end{eqnarray*}
and the following approximations of the first
\begin{eqnarray*}
 (\Delta \phi)^{m}_{i+\frac{1}{2} } = \phi_{i+1}^m -  \phi_i^m\\
 (\Delta \phi)^{m+1}_{i } = \phi_{i+\frac{1}{2}}^{m+1} -  \phi_{i-\frac{1}{2}}^{m+1}
  \end{eqnarray*}
 and the second derivatives
  \begin{eqnarray*}
 \mathcal{D} (\Delta \phi)^{m}_{i+\frac{1}{2} }&=&MM [ (\Delta \phi)^{m}_{i+\frac{3}{2}}  - (\Delta  \phi)^{m}_{i+\frac{1}{2}},  \frac{1}{2} \left((\Delta  \phi)^{m}_{i+\frac{3}{2}}  - (\Delta  \phi)^{m}_{i-\frac{1}{2}} \right),\\
 &&
 (\Delta  \phi)^{m}_{i+\frac{1}{2}}  - (\Delta  \phi)^{m}_{i-\frac{1}{2}} ]\\
 \mathcal{D} (\Delta \phi)^{m+1}_{i } &=& MM [ (\Delta \phi)^{m+1}_{i+1}  - (\Delta  \phi)^{m+1}_{i},  \frac{1}{2} \left((\Delta  \phi)^{m+1}_{i+1}  - (\Delta  \phi)^{m+1}_{i-1} \right),\\
 &&
 (\Delta  \phi)^{m+1}_{i}  - (\Delta  \phi)^{m+1}_{i-1} ]
    \end{eqnarray*}
   with the Min-Mod function
     \begin{eqnarray*}
     MM(x_1,x_2,x_3) =  \left\{
     \begin{array}{ll}
           \mbox{min}_j \{ x_j\} , &\mbox{if  all} \; x_j >0\\
          \mbox{max}_j \{ x_j\} ,&  \mbox{if  all} \;  x_j <0\\
          0, & \mbox{otherwise} .        \end{array}
          \right.
            \end{eqnarray*}
The limiter is used to deal with the possible appearance of discontinuites.

\begin{remark}
For the above second order scheme a CFL condition has to be fulfilled:
$$
\frac{\Delta t}{\Delta x}  \vert   \lambda_{max} \vert  \le \frac{1}{2}
$$
where $\lambda_{max}$ is the maximal (in absolute value) eigenvector of 
$ \frac{\partial H}{\partial \phi_x} (\phi, \phi_x )$.
Thus, for the Hamilton-Jacobi model the choice of the time step depends on the values of the gradient $\partial_x u$ and might be very small for very sharp gradients. This could be avoided by using, for example, equation (\ref{merge}).
\end{remark}

\begin{remark}
We note that  using the above described  second order method for the Aw-Rascle equations with  situations involving  contact discontinuities 
gives, among other problems, quite diffusive results. This is observed for  classical numerical methods for hyperbolic equations
as well, see \cite{CG07}. For a strategy to compute   the contact discontinuities  in a more  accurate and efficient way we refer to 
\cite{Tor,CG07}. 
\end{remark}

 \subsection{Numerical examples}

For the numerical simulations we consider the equations (\ref{macro}), (\ref{HJ}) with coefficients (\ref{a}), (\ref{b}) respectively and the constants $H=1, v_{ref}=1$, i.e. $\rho_{max} =1$.
The behavior of the solutions to the macroscopic equations is investigated in four different test scenarios. 
To illustrate the performance of the scheme described above the results are presented with two different  mesh sizes $\Delta x = 0.01$ and $\Delta x = 0.001$.
All test cases start with Riemann problems of the following form:
\begin{eqnarray*}
 \phi (x,0) =  \left\{  \begin{array}{ll}
           \phi_l , &\mbox{for} \; x <x_0\\
         \phi_r ,&  \mbox{for } \;  x>x_0
                 \end{array}
          \right.
            \end{eqnarray*}
where 
 \begin{eqnarray*}
 \phi _{l/r}=  \left(  \begin{array}{ll}
           \rho_{l/r} \\
         u_{l/r}
                 \end{array}
          \right)
            \end{eqnarray*}
are given as initial data. 

\begin{figure}
\begin{center}
\epsfig{file=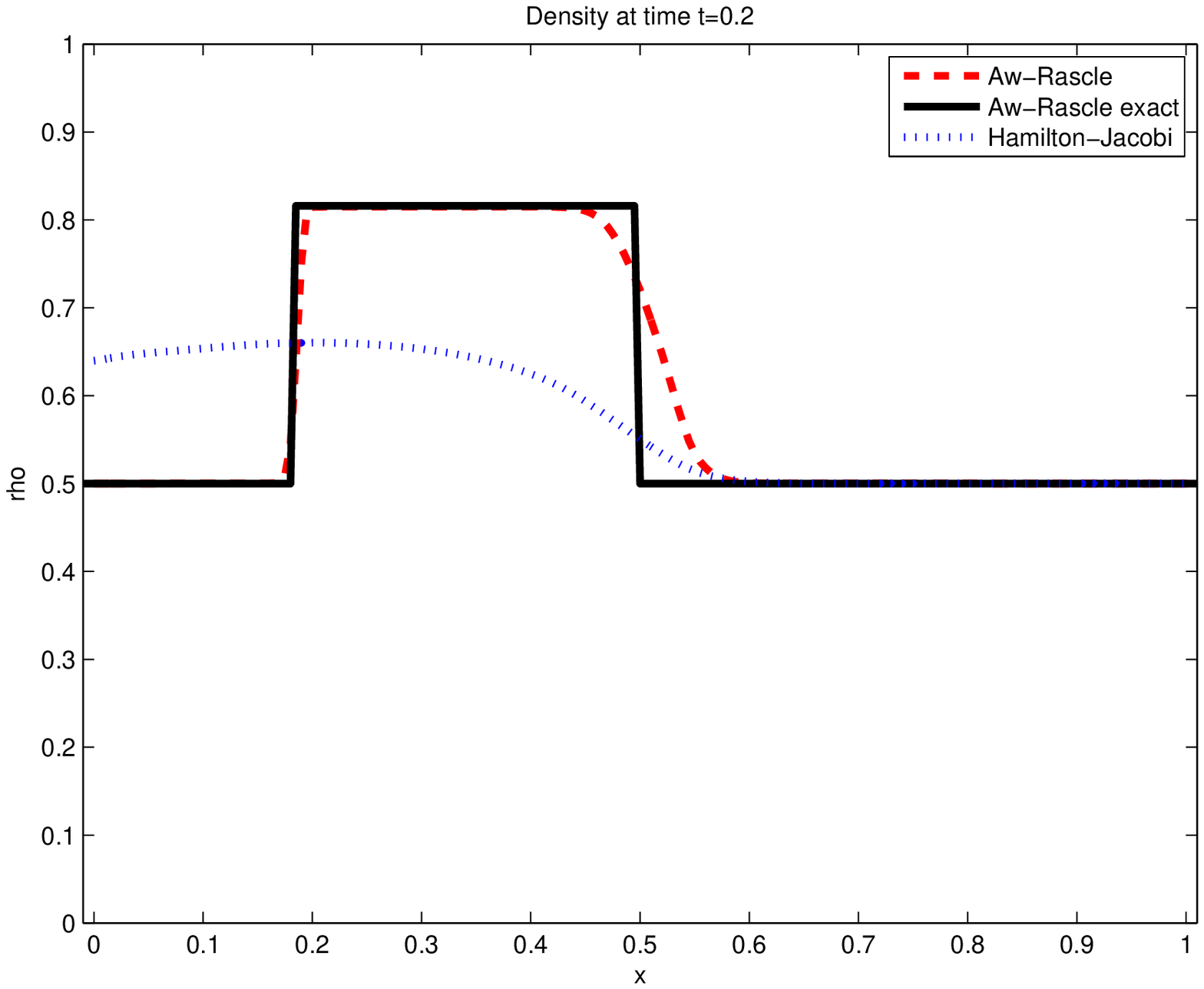,height=.25\textheight}
\epsfig{file=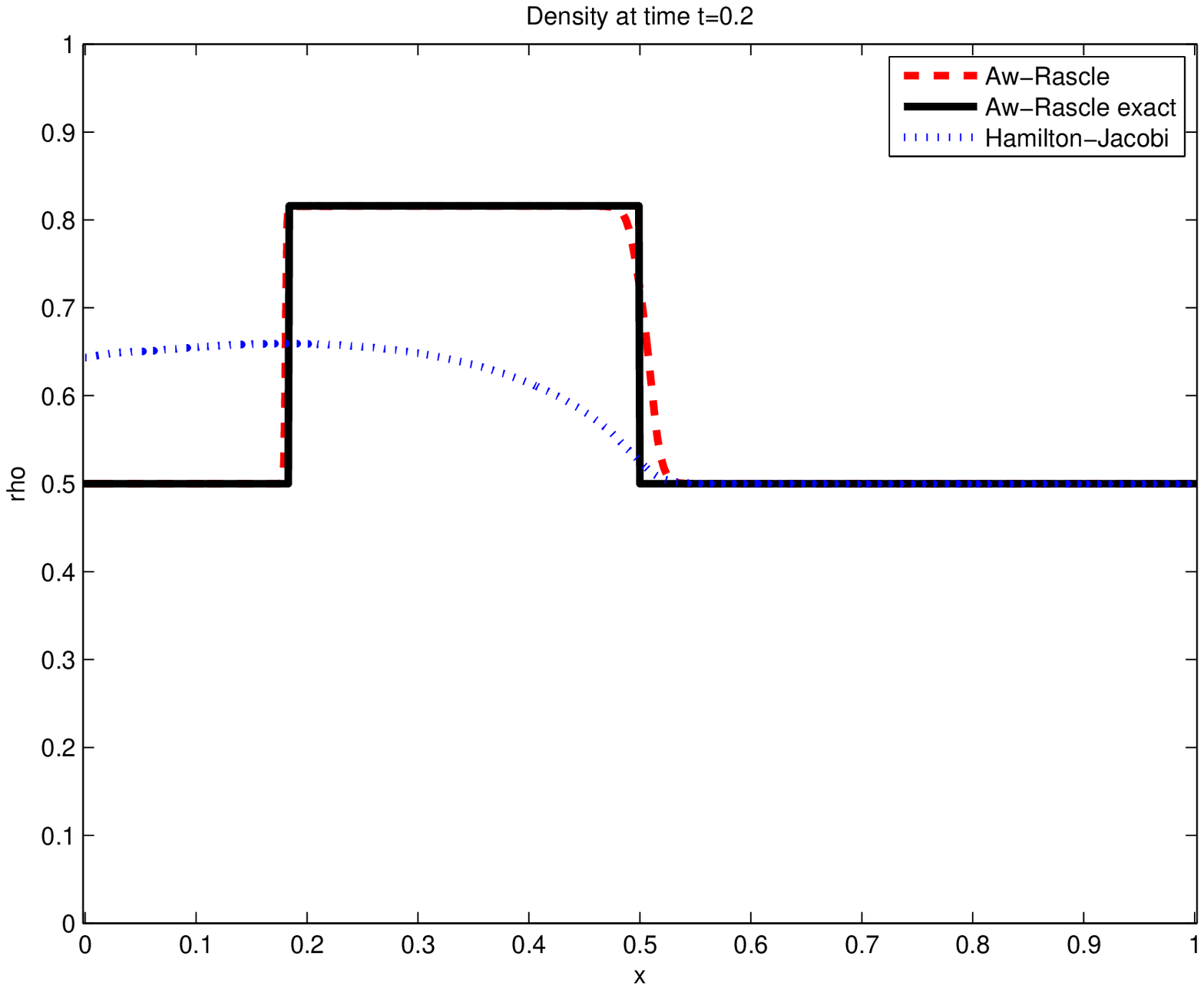,height=.25\textheight}
\end{center}
\caption{Density $\rho$ at $t=0.2$ for the Riemann problem with $ \rho _l = 0.5  , \;   u_l=1, \; 
 \rho _r =0.5 , \;   u_r= 0$ and $x_0=0.5$.}
\label{ex1}
\end{figure}

{\bf Example 1:}
In the first example the end of a traffic jam is considered. 
Thereby fast cars approach from the left a group of cars at rest on the right.
The corresponding data is given by
\begin{align*}
& \rho _l = 0.5  , \;   u_l=1,
&& 
& \rho _r =0.5 , \;   u_r= 0
 \end{align*}
 and $x_0=0.5$.
 For the Rascle model the computations are performed in conservative form using the variables $(\rho,y = \rho (u- \mbox{ln} (1-\rho))) $ to obtain the correct shock speeds.
 The numerical results are shown in Figure \ref{ex1}.
 The exact solution of the Rascle model (solid line) is given by a shock-wave moving to the right followed by a stationary contact-discontinuity, see \cite{AR98}.
The numerical results for the Hamilton-Jacobi model (dotted line) show a faster braking of the approaching cars.
This leads to a faster back-traveling wave and a less dense congested state. 
About the numerical aspects, the diffusion at the contact discontinuity is reduced by the finer grid, whereas the resolution of the shock in the Rascle model (dashed line) remains satisfactory.

\begin{figure}
\begin{center}
\epsfig{file=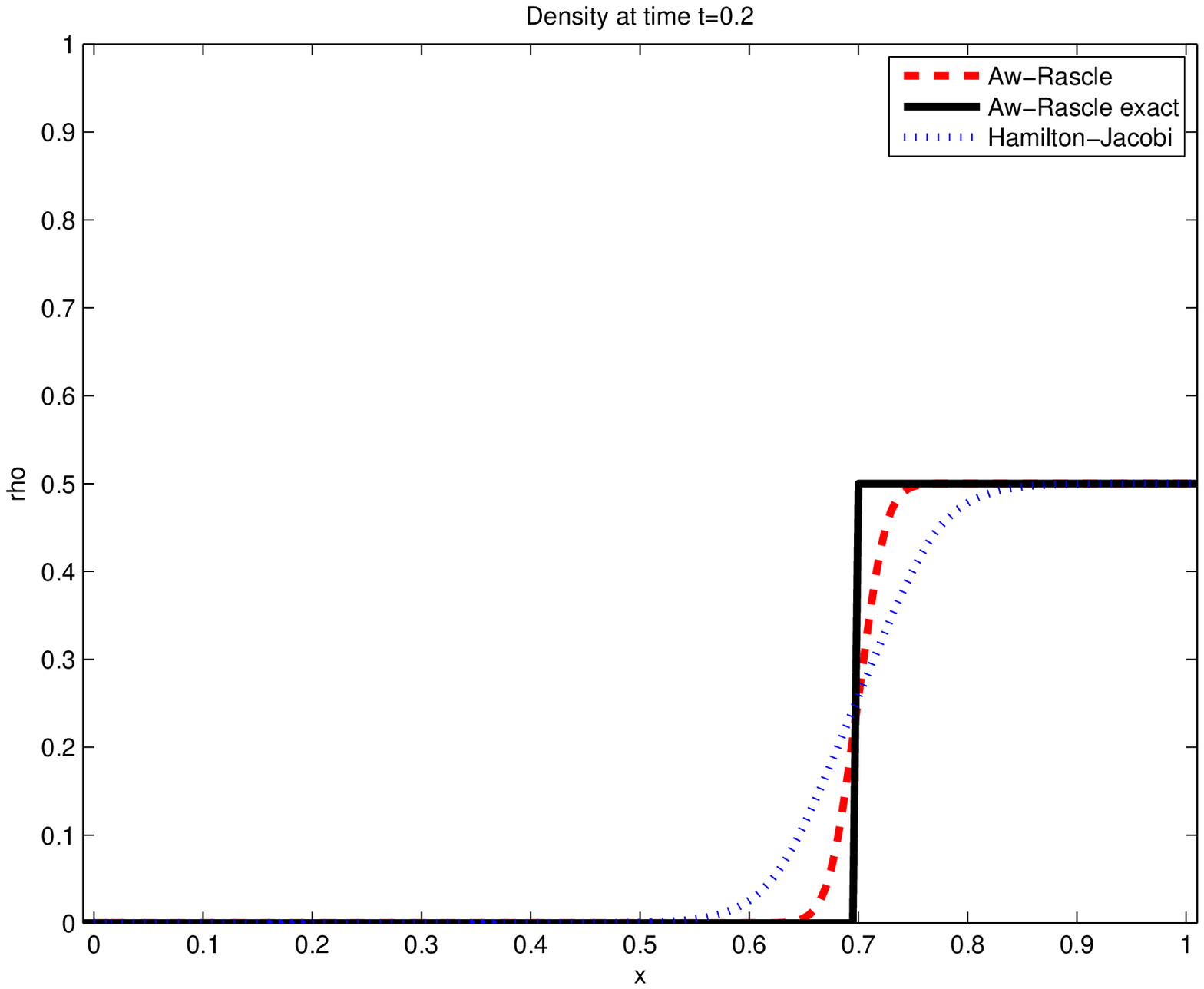,height=.25\textheight}
\epsfig{file=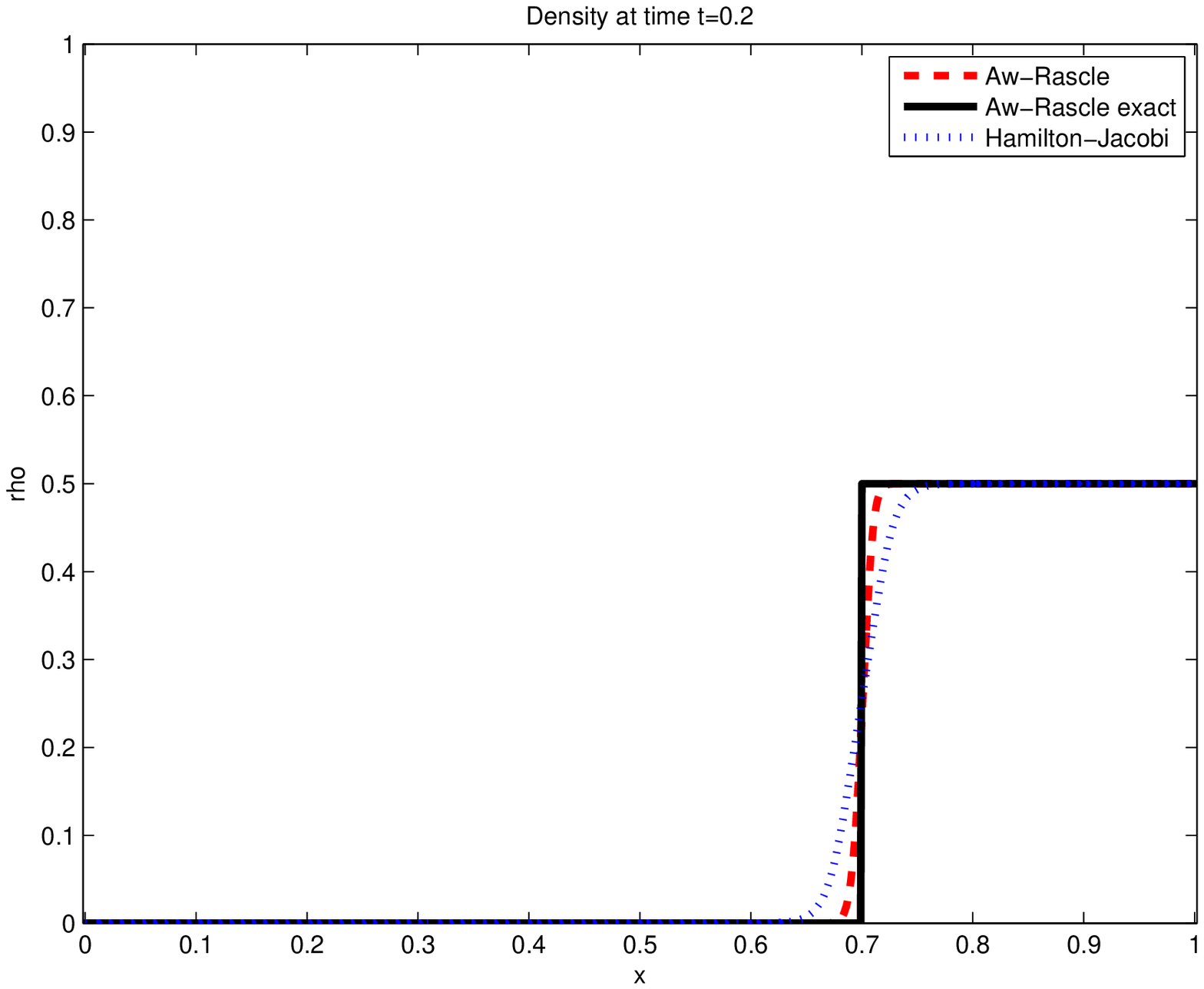,height=.25\textheight}
\end{center}
\caption{Density $\rho$ at $t=0.2$ for the Riemann problem with $ \rho _l = 0  , \;   u_l=1, \; 
 \rho _r =0.5 , \;   u_r= 1$ and $x_0=0.5$.}
\label{ex2}
\end{figure}

{\bf Example 2:}
Now the tail of a group of moving cars followed by an empty road is studied. The initial states are chosen as
\begin{align*}
 &\rho _l =  0 , \;   u_l= 1
 &&\text{and}
 &\rho _r =  0.5, \;   u_r=1\ ,
 \end{align*}
 with the discontinuity at $x_0=0.5$.
   As shown in Figure \ref{ex2}, the exact solution of the Rascle model (solid line) is given by a single contact-discontinuity moving at the speed of the leading cars.
   This behavior is captured well by the numerical scheme (dashed line) and holds also true for the Hamilton-Jacobi model. 
   In both cases the cars are not influenced by the free space behind them and are thus following the constant state in front.

\begin{figure}
\begin{center}
\epsfig{file=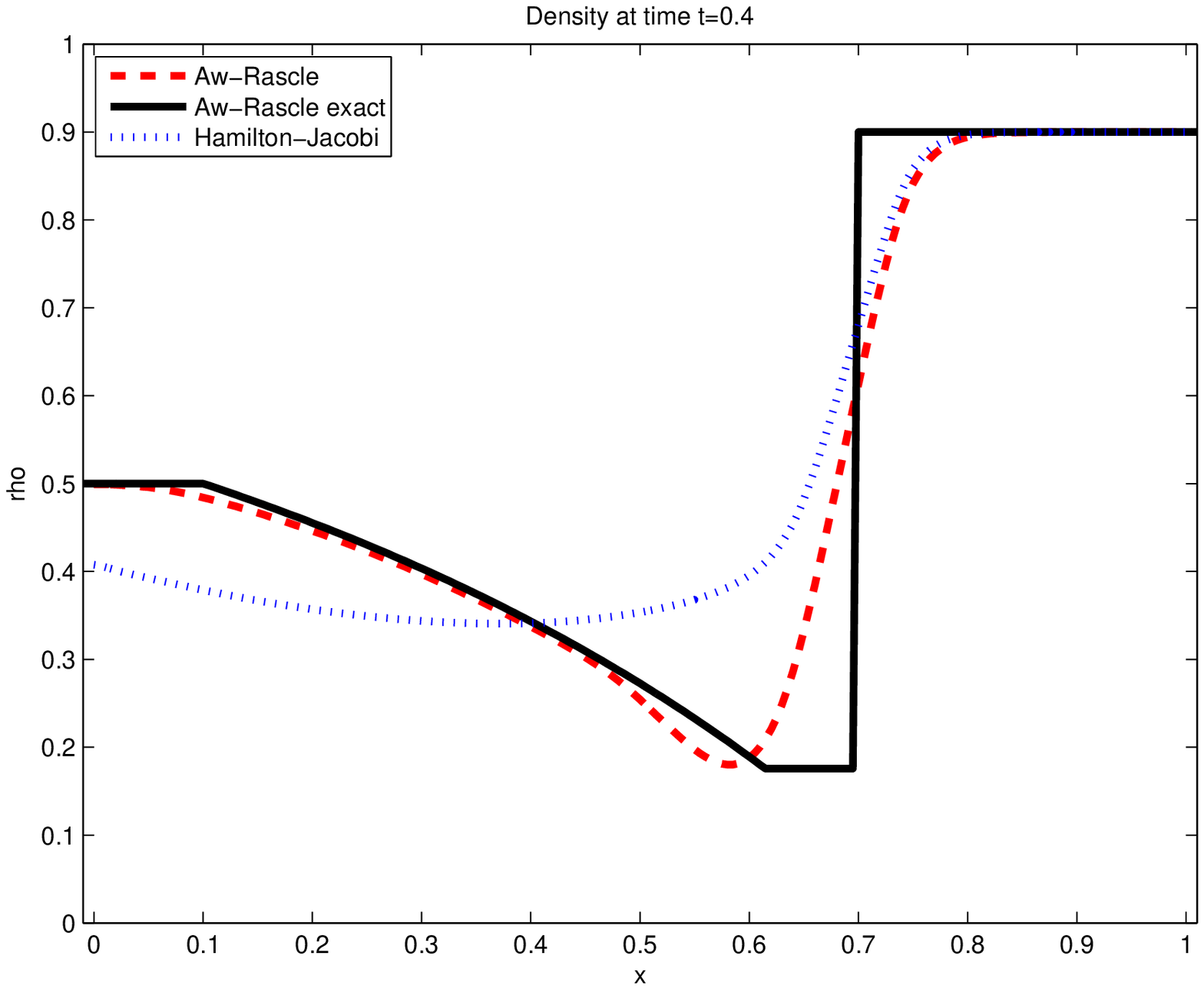,height=.25\textheight}
\epsfig{file=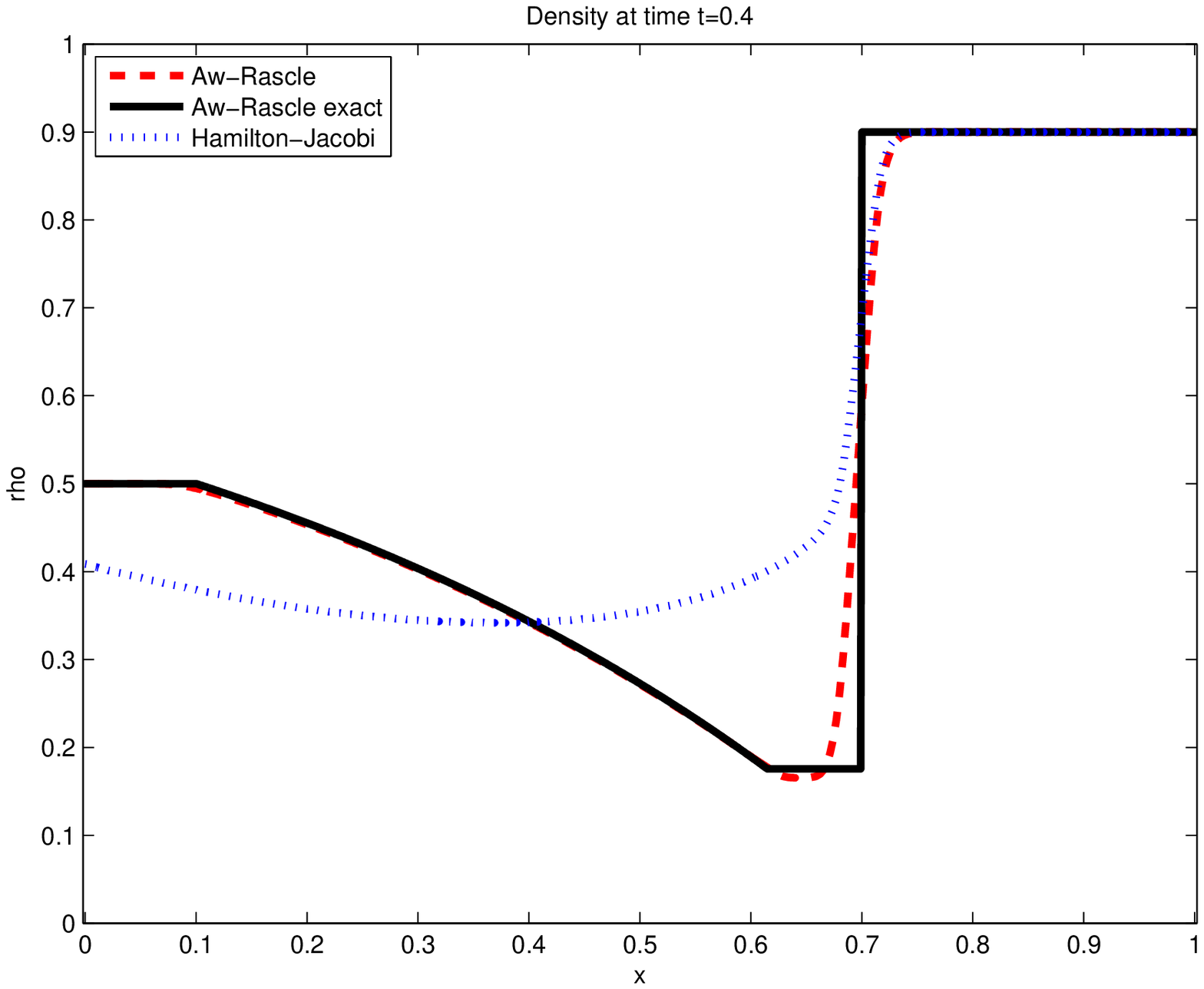,height=.25\textheight}
\end{center}
\caption{Density $\rho$ at $t=0.4$ for the Riemann problem with $ \rho _l = 0.5  , \;   u_l=0, \; 
 \rho _r =0.9 , \;   u_r= 0.5$ and $x_0=0.5$.}
\label{ex3}
\end{figure}

{\bf Example 3:}
Here we consider a group of faster vehicles escaping from slower ones in behind.
Therefore we chose
\begin{align*}
 &\rho _l = 0.5 , \;   u_l = 0
&&\text{and} 
 &\rho _r = 0.9 , \;   u_r =0.5
 \end{align*}
 on the left and right of $x_0=0.5$.
 In Figure \ref{ex3} the corresponding solutions are plotted.
 The exact solution of the Rascle model (solid line) consists of a left going rarefaction wave  and a contact-discontinuity moving to the right. 
 As the drivers of the Hamilton Jacobi model (dotted line) tend to accelerate faster than those of the Rascle model (dashed line), the arising gap is less distinct.
 Thus a more homogeneous state is reached on the left.
 By increasing the number of grid points only the resolution of the contact discontinuity is improved.

\begin{figure}
\begin{center}
\epsfig{file=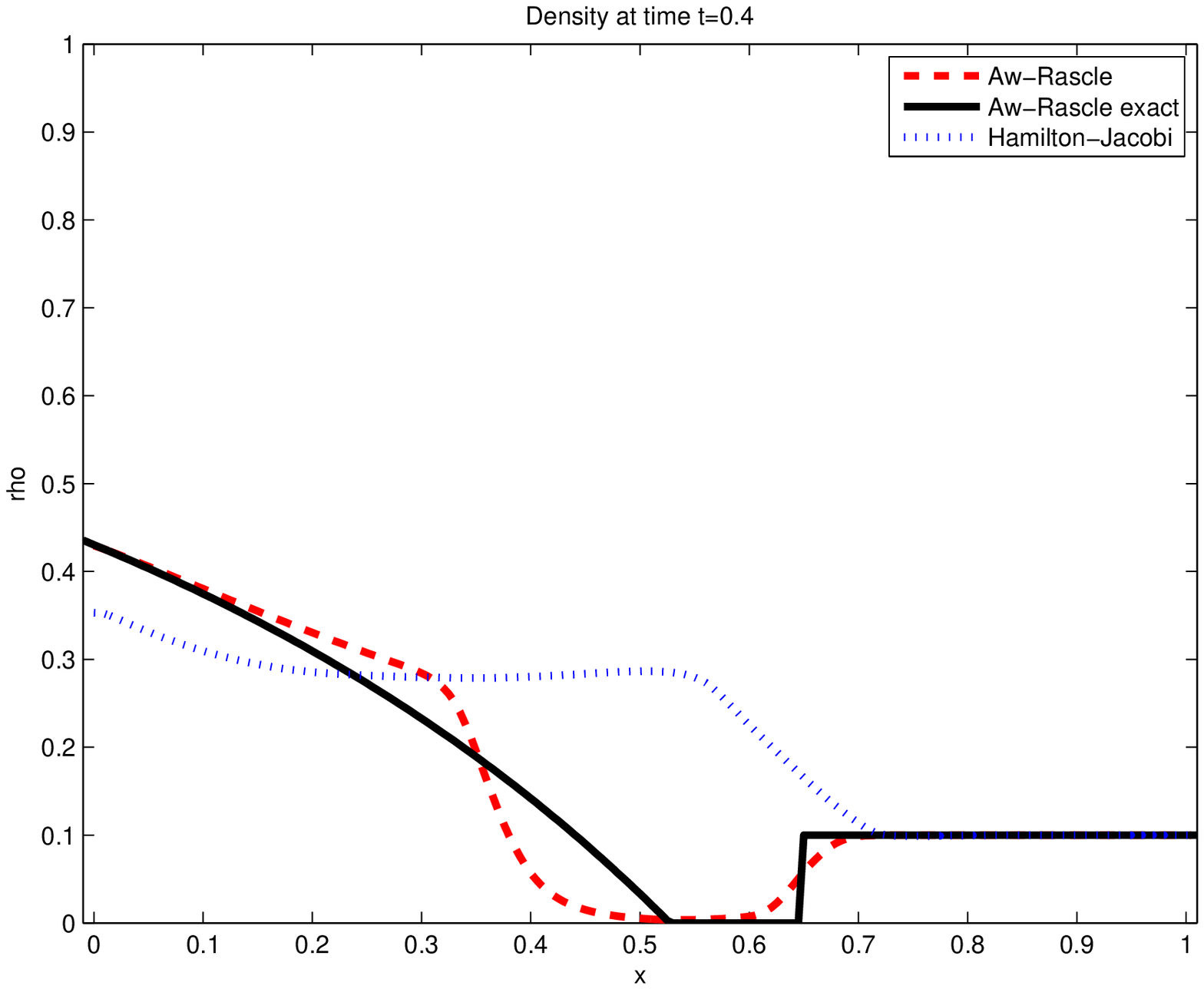,height=.25\textheight}
\epsfig{file=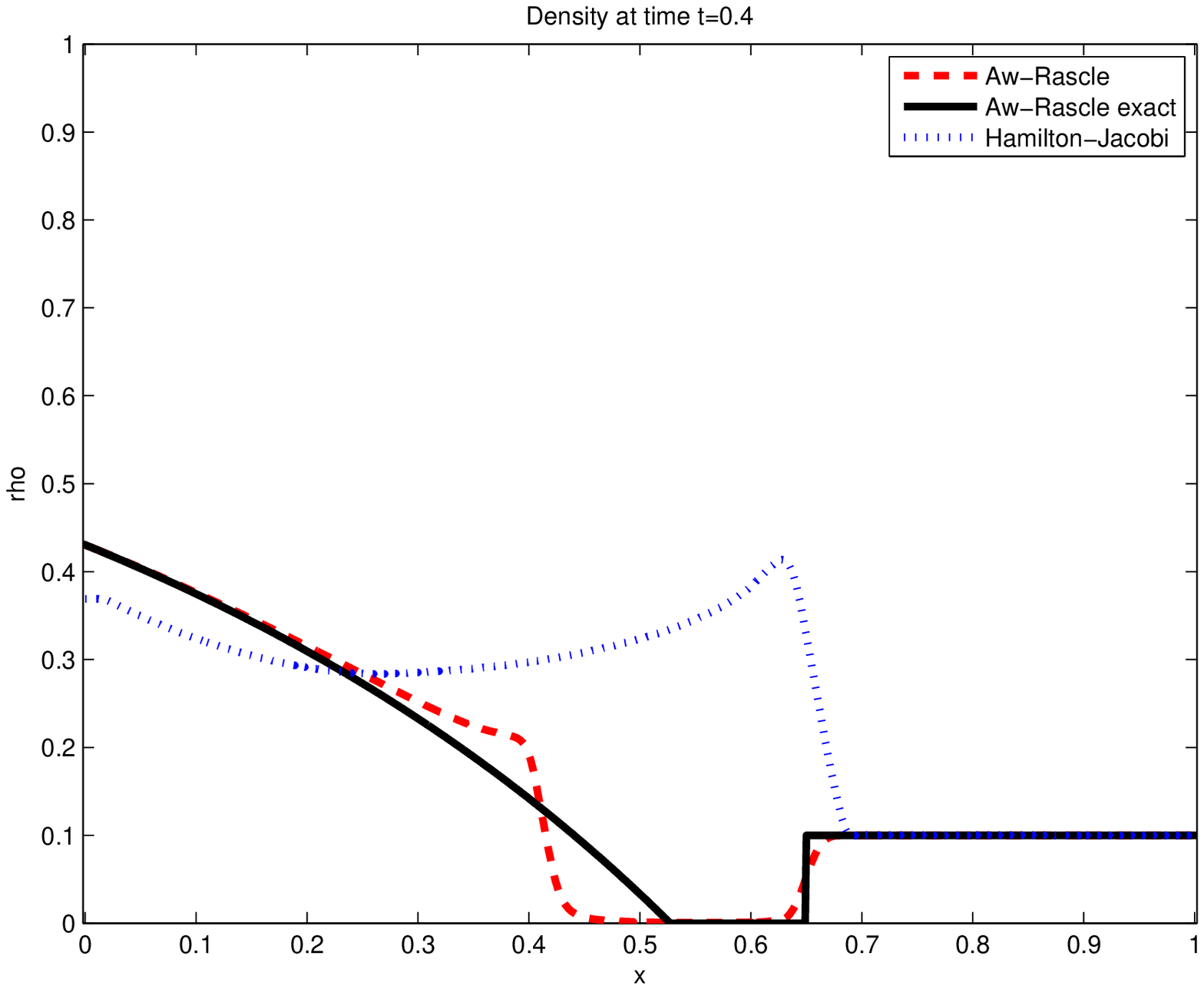,height=.25\textheight}
\end{center}
\caption{Density $\rho$ at $t=0.5$ for the Riemann problem with $ \rho _l = 0.5  , \;   u_l=0, \; 
 \rho _r =0.1 , \;   u_r= 1$ and $x_0=0.25$.}
\label{ex4}
\end{figure}

{\bf Example 4:}
Finally we consider an example similar to the above one, but now with faster cars on the right.
The data is given as
\begin{align*}
& \rho _l = 0.5 , \;   u_l= 0
,&&
& \rho _r = 0.1 , \;   u_r=1
 \end{align*}
  and $x_0=0.25$. 
  The exact solution of the Rascle model (solid line, Figure \ref{ex4}) is given by a rarefaction wave connected to a vacuum state, which is followed by a contact-discontinuity moving to the right.
  Here a difference to the numerical solution (dashed line) is observed. The applied scheme fails to properly capture the fake wave connecting the rarefaction wave to the vacuum state.
  The artificial jump can not be reduced by an increase of the computational accuracy.
  In the Hamilton Jacobi model (dotted line) no such vacuum state arises, since the drivers tend to accelerate faster. 

In the above examples the wave fronts for the Hamilton-Jacobi model are smeared compared to the Aw-Rascle model as expected.
In particular, Example 1 shows a stronger breaking for the Hamilton-Jacobi model and example 3 shows a faster acceleration of the vehicles keeping contact with the leading cars.

\begin{rem}
One also observes comparing the coarse and fine grid numerical solution, that the Hamilton-Jacobi equations are already well approximated by the coarse grid solution. Only example 4 shows a further steepening of the solution by refining the mesh. In general, the Rascle type conservation law is well approximated by the scheme except some smearing of the contact discontinuities.
The only exception is the vacuum wave in example 4, where a non-physical jump is generated. 
Numerical difficulties at vacuum states are discussed e.g.  in \cite{Toro}.
\end {rem}

\begin{remark}
The numerical solution of the  hyperbolic Aw-Rascle model is sensitive to the choice of variables. 
Example 1 (a solution with a shock) is computed using conservative variables $(\rho, y= \rho(u - \mbox{ln} (1-\rho))$ to ensure  the correct intermediate state. 
Example 2,3,4 have been computed in $(\rho, \rho u)$ variables, since no shocks appear.
Although this choice of variables improves the resolution of the contact discontinuity it remains rather diffusive.
As mentioned above using the methods described in \cite{Tor,CG07} a sharp resolution of the contact discontinuities can be obtained.
Nevertheless, we plotted in the above figures for comparison the solutions using the scheme described in Section 5.1.
\end {remark}

\subsubsection*{Conclusions}

\begin{itemize}
\item
The paper contains the derivation of  two classes of macroscopic models from kinetic equations.
The type  of equation one obtains does not depend on the fact whether an integro-differential equation or a Fokker-Planck type model is used, but rather on the fact 
which interaction rule is chosen.
\item
In certain cases a Hamilton-Jacobi term can be derived in the momentum equations instead of the classical Rascle term.
\item
Numerical investigation using a suitable second order method have been used to investigate the behavior of the solutions showing a smearing effect 
of the wave fronts for the Hamilton-Jacobi equations.
\item
Further investigations will include the derivation of suitable relaxation  terms from kinetic models and multiphase traffic equations.
\end{itemize}

\subsubsection*{Acknowledgments}

The present work has been supported by DFG KL  1105/16-1 and the DAAD PhD Program  MIC.


\begin{thebibliography}{10}

\bibitem{aw00}
{\sc A.~Aw}, {\em Mod\`eles hyperboliques de trafic automobile}, PhD thesis,
  Nice, 2001.
  
\bibitem{AKMR03}
A. Aw, A. Klar, T. Materne, and M. Rascle,
{\em Derivation of Continuum Traffic Flow Models from Microscopic
Follow-the-Leader Models.}
SIAM J. Appl. Math. 63/1, 259-278, 2002
\bibitem{AR98}
{\sc A.~Aw and M.~Rascle}, {\em Resurrection of second order models of traffic
  flow?}, SIAM J. Appl. Math., 60 (2000), pp.~916--938.
  
  
 \bibitem{BL03}
{\sc S. Bryson, D. Levy}, {\em Central schemes for multidimensional Hamilton-Jacobi equations}
SIAM Sci. Comp. 25, 3, 767-791, 2003

  \bibitem{CDP}
Carrillo, J.A., D'Orsogna, M.R., Panferov, V.: \emph{Double milling in
self-propelled swarms from kinetic theory}. Kinetic and Related
Models, \textbf{2} (2009), pp. 363-378.



 \bibitem{Deg}
{\sc F. Berthelin, P. Degond, M. Delitla, M. Rascle}, {\em A model for the formation and evolution of traffic jams}
Arch. Rat. Mech. Anal. 187, 185-220, 2008

\bibitem{Ber}
{\sc F. Berthelin, P. Degond, V. Le Blanc, S. Moutari, J. Royer, M. Rascle},  
{\em A Traffic-Flow Model with Constraints for the Modeling of Traffic Jams}, Mathematical Models and Methods in Applied Sciences  18,   1269-1298, 2008

\bibitem{CG07}
{\sc C. Chalons, P. Goatin},  
{\em Transport-equilibrium schemes for computing contact discontinuities in traffic flow modelling}, 
Commun. Math. Sci. Volume 5,3, 533-551, 2007

\bibitem{GKMW03}
{\sc M.~G\"unther, A.~Klar, T.~Materne, and R.~Wegener},
\newblock {\em Multivalued fundamental diagrams and stop and go waves for continuum
  traffic flow equations}.
\newblock SIAM J. Appl. Math. 64/2, 468-483, 2003

 
\bibitem{Gre00}
{\sc J.~Greenberg}, {\em Extension and amplification of the {A}w-{R}ascle
  model}, SIAM J. Appl. Math., 62 (2001), pp.~729--745.

\bibitem{Hel95B}
{\sc D.~Helbing}, {\em Gas-kinetic derivation of {N}avier-{S}tokes-like traffic
  equation}, Physical Review E, 53 (1996), pp.~2366--2381.

\bibitem{IH}
{\sc M. Herty, R. Illner},  {\em On stop and go waves in dense traffic}, 
Kinetic and Related Models (KRM),1, 2008, 437-452


\bibitem{IKM03}
{\sc R.~Illner, A.~Klar, and T.~Materne}, {\em Vlasov-fokker-planck models for
  multilane traffic flow}, Comm. Math. Sci., 1 (2003), pp.~1--12.

\bibitem{KW97}
{\sc A.~Klar and R.~Wegener}, {\em Enskog-like kinetic models for vehicular
  traffic}, J. Stat. Phys., 87 (1997), pp.~91--114.

\bibitem{KW981}
{\sc A.~Klar and R.~Wegener},  {\em A hierachy of models
  for multilane vehicular traffic {I}: Modeling}, SIAM J. Appl. Math., 59
  (1998), pp.~983--1001.

\bibitem{KW00}
 {\sc A.~Klar and R.~Wegener}, {\em Kinetic derivation
  of macroscopic anticipation models for vehicular traffic}, SIAM J. Appl.
  Math., 60 (2000), pp.~1749--1766.

\bibitem{Nel95}
{\sc P.~Nelson}, {\em A kinetic model of vehicular traffic and its associated
  bimodal equilibrium solutions}, Transport Theory and Statistical Physics, 24
  (1995), pp.~383--408.

\bibitem{PF75}
{\sc S.~Paveri-Fontana}, {\em On {B}oltzmann like treatments for traffic flow},
  Transportation Research, 9 (1975), pp.~225--235.

\bibitem{Pay79}
{\sc H.~Payne}, {\em {FREFLO}: A macroscopic simulation model of freeway
  traffic}, Transportation Research Record, 722 (1979), pp.~68--75.

\bibitem{PH71}
{\sc I.~Prigogine and R.~Herman}, {\em Kinetic Theory of Vehicular Traffic},
  American Elsevier Publishing Co., New York, 1971.

\bibitem{Tor}
{\sc E.F. Toro}, {\em Riemann solvers and numerical methods for fluid dynamics},
 Springer, Berlin, Heidelberg, 2009.


\bibitem{Toro}
{\sc E.F. Toro},{\em Shock-Capturing Methods for Free-Surface Shallow Flows},
John Wiley, 2001

\bibitem{Whi74}
{\sc G.~Whitham}, {\em Linear and Nonlinear Waves}, Wiley, New York, 1974.

\bibitem{WL}
{\sc W. Leutzbach, R. Wiedemann}, {\em  Development and application of traffic simulation models at Karlsruhe Institut fuer Verkehrswesen},
Traffic engineering and control, May, 20-278, 1986
\end{thebibliography}
\end{document}